\documentclass[a4paper,twoside]{article}
\usepackage{amssymb,amsmath,amscd,amsthm}

\input{epsf}

\newcommand{\str}{^\prime}
\newcommand{\res}[1]{\!\mid_{#1}}
\newcommand{\card}[1]{{\mid\! #1 \!\mid}}
\newcommand{\pro}[2]{\langle #1, #2 \rangle}

\newcommand{\norm}[1]{\| #1 \|}

\def\lolra{\Longleftrightarrow}
\def\lora{\Longrightarrow}
\def\lra{{\,\Leftrightarrow\,}}
\def\ra{{\,\Rightarrow\,}}

\def\C{{\mathbb C}}
\def\N{{\mathbb N}}
\def\Nspos{{\N_{>0}}}
\def\P{{\mathbb P}}
\def\Q{{\mathbb Q}}
\def\R{{\mathbb R}}
\def\Z{{\mathbb Z}}
\def\Rpos{{\R_{\geq 0}}}

\def\F{{\cal F}}

\def\M{{\cal M}}
\def\calN{{\cal N}}
\def\J{{\cal R}}
\def\S{{\cal S}}
\def\TT{{\mathbb T}}

\def\calU{{\cal U}}
\def\V{{{\cal V}}}

\def\conv{{\rm conv}}
\def\aff{{\rm aff}}
\def\lin{{\rm lin}}
\def\pos{{\rm pos}}
\def\rvol{{\rm rvol}}

\def\vol{{\rm vol}}
\def\codim{{\rm codim}}
\def\relint{{\rm relint}}

\def\Aut{{\rm Aut}}
\def\Hom{{\textup{Hom}}}

\def\Cl{\textup{Cl}}

\def\NR{{N_{\R}}}
\def\NQ{{N_{\Q}}}
\def\MR{{M_{\R}}}
\def\MQ{{M_{\Q}}}

\def\fan{\triangle}

\def\taus{{\tau\str}}

\def\rand{\partial}
\def\randp{{\rand P}}

\def\rays{{\fan(1)}}

\def\intr{{\rm int}}

\def\randPM{{\rand P \cap M}}
\def\intP{\intr P}

\def\Vt{{\V_\tau}}

\def\vt{{v_\tau}}

\def\vts{{v_\taus}}

\def\Zt{{\cal Z}}

\def\fb{{\F(\calU)}}
\def\GFt{\Z/3\Z}

\newtheorem*{theorem*}{Theorem}
\newtheorem*{corollary*}{Corollary}
\newtheorem{theorem}{Theorem}[section]
\newtheorem{definition}[theorem]{Definition}
\newtheorem{remark}[theorem]{Remark}
\newtheorem{example}[theorem]{Example}
\newtheorem{corollary}[theorem]{Corollary}
\newtheorem{proposition}[theorem]{Proposition}
\newtheorem{lemma}[theorem]{Lemma}

\newtheorem*{acknowledgement}{Acknowledgement}

\title{Complete toric varieties with\\reductive automorphism group}

\author{{\sc Benjamin Nill} \\
\small  {\em Mathematisches Institut, Universit\"at T\"ubingen}   \\
\small  {\em Auf der Morgenstelle 10,  72076  T\"ubingen, Germany}  \\
\small  {\em e-mail: nill@everest.mathematik.uni-tuebingen.de} \\
 }

\begin{document}

\date{}

\maketitle

\begin{abstract}
We give equivalent and sufficient criteria for the automorphism group 
of a complete toric variety, respectively a Gorenstein toric Fano variety, to be reductive. 
In particular we show that the automorphism group of a Gorenstein toric Fano variety is reductive, if 
the barycenter of the associated reflexive polytope is zero. Furthermore 
a sharp bound on the dimension of the reductive automorphism group of a complete toric variety 
is proven by studying the set of Demazure roots. 
\end{abstract}

\section*{Introduction}
\label{intro}

There is an important obstruction to the existence of an Einstein-K\"ahler metric on a nonsingular Fano 
variety:

\begin{theorem*}[Matsushima 1957]
If a nonsingular Fano variety $X$ admits an Einstein-K\"ahler metric, then $\Aut(X)$ is a reductive algebraic group.
\end{theorem*}

In 1983 Futaki introduced the so called {\em Futaki character}, 
whose vanishing is another important necessary condition for the existence of an Einstein-K\"ahler metric. 
For a nonsingular toric Fano variety with reductive automorphism group there is an explicit criterion 
(see \cite[Cor. 5.5]{Mab87}):

\begin{theorem*}[Mabuchi 1987]
Let $X$ be a nonsingular toric Fano variety with $\Aut(X)$ reductive.

The Futaki character of $X$ vanishes if and only if 
the barycenter of $P$ is zero, where $P$ is the associated reflexive polytope, i.e., 
the fan of normals of $P$ is associated to $X$. 
\end{theorem*}

In \cite[Thm. 1.1]{BS99} Batyrev and Selivanova were able to give a sufficient criterion for the 
existence of an Einstein-K\"ahler metric.

\begin{theorem*}[Batyrev/Selivanova 1999]
Let $X$ be a nonsingular toric Fano variety. We denote by $P$ the associated reflexive polytope. 

If $X$ is symmetric, i.e., the group of lattice automorphisms leaving $P$ invariant 
has no non-zero fixpoints, then $X$ admits an Einstein-K\"ahler metric.
\end{theorem*}

In particular they got as 
a corollary \cite[Cor. 1.2]{BS99} that the automorphism group of such a symmetric $X$ is reductive. 
Expressed in combinatorial terms this just means 
that the set of lattice points in the relative interiors of facets of $P$ is centrally symmetric. So 
they asked whether a direct proof for this result exists. 
Indeed there is even the following generalisation with a simple combinatorial proof (see Theorem \ref{main}(1)):

\begin{theorem*}
Let $X$ be a complete toric variety.

If the group of automorphisms of the associated fan has no non-zero fixpoints, then $\Aut(X)$ is reductive.
\end{theorem*}

Motivated by above results it was conjectured by Batyrev that in the case of a nonsingular 
toric Fano variety already the vanishing of the barycenter of the associated reflexive polytope were sufficient 
for the automorphism group to be reductive. Indeed there is even the following more general result that has 
a purely convex-geometrical proof (see Theorem \ref{main}(2{\rm i})):

\begin{theorem*}
Let $X$ be a Gorenstein toric Fano variety. 

If the barycenter of the associated reflexive polytope is zero, then $\Aut(X)$ is reductive.
\end{theorem*}

Only very recently Xu-Jia Wang and Xiaohua Zhu could prove that the vanishing of the Futaki character 
is even sufficient for the existence of an Einstein-K\"ahler metric in the toric case (see \cite[Cor. 1.3]{WZ03}):

\begin{theorem*}[Wang/Zhu 2003]
Let $X$ be a nonsingular toric Fano variety. 

Then $X$ admits an Einstein-K\"ahler metric if and only the Futaki character of $X$ vanishes.
\end{theorem*}

Combining the previous results this yields a generalisation of the above theorem of Mabuchi that 
is also implicit in \cite[Lemma 2.2]{WZ03}:

\begin{corollary*}
Let $X$ be a nonsingular toric Fano variety.

Then $X$ admits an Einstein-K\"ahler metric if and only if 
the barycenter of $P$ is zero, where $P$ is the associated reflexive polytope. 
\end{corollary*}

It is now conjectured by Batyrev 
that this result may also hold in the singular case of a Gorenstein toric Fano variety.

The paper is organised as follows:

In the first section the notation is fixed and basic definitions are given.

The second section deals with the automorphism group of a $d$-dimen\-sional complete toric variety $X$. Here 
the set of roots $\J$ plays a crucial part in determining the dimension of the identity component and 
whether the group is reductive (see Prop. \ref{autdim}). 
Using results of Cox in \cite{Cox95} we construct pairwise orthogonal families of roots, 
so called $\S$-root bases, that parametrize the set of semisimple roots $\S$ in a convenient way, 
where $\S := \J \cap -\J$. In Prop. \ref{charac} it is shown that 
$X$ is isomorphic to a product of projective spaces if and only if there are $d$ linearly independent semisimple roots. 
When $\Aut(X)$ is reductive, $\dim \Aut^\circ(X) > d^2 - 2$ is a sufficient condition for that 
(see Theorem \ref{root-dim}). 

In the third section we more closely examine the case of a 
$d$-dimensional Goren\-stein toric Fano variety $X$ associated to a reflexive polytope $P$ (see \cite{Nil04}). 
Here a root of $X$ is just a lattice point in 
the relative interior of a facet of $P$, so the results of the previous section have a direct geometric 
interpretation. In particular we obtain that $P$ has at most $2 d$ facets containing roots of $P$, 
with equality if and only if $X$ is the product of $d$ projective lines (see Corollary \ref{wuerfel}). Furthermore 
the intersection of $P$ with the space spanned by all semisimple roots 
is a reflexive polytope associated to the product of projective spaces (see Theorem \ref{intersec}). 

In the fourth section we present equivalent and sufficient criteria for the automorphism group 
of a complete toric variety, respectively a Gorenstein toric Fano variety, 
to be reductive (see Theorem \ref{main}).

In the last section we are concerned with $d$-dimensional centrally symmetric reflexive polytopes. 
In particular we finish in Theorem \ref{rootcentral}(3) the proof of 
\cite[Thm. 6.4]{Nil04} saying that such a polytope has 
at most $3^d$ lattice points, with equality if and only if the associated toric variety is the product of 
$d$ projective lines.

\begin{acknowledgement}{\rm
The author would like to thank his thesis advisor 
Professor Victor Batyrev for posing problems, his advice 
and encouragement. The author would also like to thank Professor Maximillian Kreuzer for 
the support with the computer package PALP, the classification data and many examples.

The author was supported by DFG, Forschungsschwerpunkt "Globale Methoden in der komplexen Geometrie". 
This work is part of the author's thesis.}
\end{acknowledgement}

\section{Notation and basic definitions}

In this section we shortly repeat the standard notation for polytopes and toric varieties, as 
it can be found in \cite{Ewa96}, \cite{Ful93} or \cite{Oda88}. In \cite{Bat94} reflexive polytopes were introduced. 

Let $N \cong \Z^d$ be a $d$-dimensional lattice and $M = \Hom_\Z(N,\Z) \cong \Z^d$ the dual lattice with 
$\pro{\cdot}{\cdot}$ the nondegenerate symmetric pairing. 
As usual, $\NQ = N \otimes_\Z \Q \cong \Q^d$ and $\MQ = M\otimes_\Z \Q \cong \Q^d$ 
(respectively $\NR$ and $\MR$) will denote the rational (respectively real) scalar extensions. 

For a subset $S$ of a real vector space let $\lin(S)$ (respectively $\aff(S)$, $\conv(S)$, $\pos(S)$) be the 
linear (respectively affine, convex, positive) hull of $S$. 
A subset $P \subseteq \MR$ is called a polytope, if it is the convex hull of finitely many points in $\MR$. 
The boundary of $P$ is denoted by $\randp$, the relative interior of $P$ by $\relint P$. When 
$P$ is full-dimensional, its relative interior is also denoted by $\intr P$. A face $F$ of $P$ is denoted by 
$F \leq P$, the vertices of $P$ form the set $\V(P)$, the facets of $P$ the set $\F(P)$. 
$P$ is called a lattice polytope, respectively rational polytope, 
if $\V(P) \subseteq M$, respectively $\V(P) \subseteq \MQ$. An isomorphism of 
lattice polytopes is an isomorphism of the associated lattices such that the induced real linear isomorphism maps 
the polytopes onto each other. 

We usually denote by $\fan$ a complete fan in $\NR$. 
The $k$-dimensional cones of $\fan$ form a set $\fan(k)$. The elements in $\fan(1)$ are called rays, 
and given $\tau \in \rays$, we let $\vt$ denote the unique generator of $N \cap \tau$.

Let $P \subseteq \MR$ be a rational $d$-dimensional polytope with $0 \in \intP$. 
Then we have the important notion of the dual polytope
\[P^* := \{y \in \NR \,:\,\, \pro{x}{y} \geq -1 \;\forall\, x \in P\},\]
that is also a rational $d$-dimensional polytope with $0 \in \intr P^*$. 
The fan $\calN_P := \{\pos(F) \,:\, F \leq P^*\}$ is called the normal fan of $P$.

Duality means $(P^*)^* = P$. 
There is a natural combinatorial correspondence between $i$-dimensional faces of $P$ and 
$(d-1-i)$-dimensional faces of $P^*$ that reverses inclusion.
 
For a facet $F \leq P$ we let $\eta_F \in \NQ$ denote the uniquely determined 
inner normal satisfying $\pro{\eta_F}{F} = -1$. We have
$$\V(P^*) = \{\eta_F \,:\, F \in \F(P)\}.$$

The dual of the product of $d_i$-dimensional polytopes $P_i \subseteq \R^{d_i}$ with $0 \in \intr P_i$ for $i=1,2$ is 
given by
\begin{equation}(P_1 \times P_2)^* = \conv(P_1^* \times \{0\}, \{0\} \times P_2^*) \subseteq \R^{d_1} \times \R^{d_2}.
\label{product}
\end{equation}

By a well-known construction a fan $\fan$ in $\NR$ 
defines a toric variety $X := X(N, \fan)$, i.e., a  normal irreducible algebraic variety over $\C$ such that 
an open embedded algebraic torus $\TT = (\C^*)^d$ acts on $X$ in extension of its own action. 

Let $P \subseteq \MR$ be a rational polytope. We define the associated toric variety 
\[X_P := X(N, \calN_P).\]

For $d$-dimensional rational polytopes 
$P_1, P_2$ equation (\ref{product}) implies
$$X_{P_1} \times X_{P_2} \cong X_{P_1 \times P_2}.$$

\begin{definition}{\rm
A complex variety $X$ is called {\em Gorenstein Fano variety}, if $X$ is projective, normal and 
its anticanonical divisor is an ample Cartier divisor.
}\end{definition}

In the toric case there is the following definition (see \cite{Bat94}):

\begin{definition}{\rm
A $d$-dimensional polytope $P \subseteq \MR$ with $0 \in \intP$ is called {\em reflexive polytope}, 
if $P$ is a lattice polytope and $P^*$ is a lattice polytope.
}\end{definition}

Especially reflexive polytopes always appear in dual pairs. 
There is the following fundamental result (see \cite{Bat94} or \cite{Nil04}):

\begin{theorem}
Under the map $P \mapsto X_P$ reflexive polytopes correspond uniquely up to isomorphism to 
Gorenstein toric Fano varieties. There are only 
finitely many isomorphism types of $d$-dimensional reflexive polytopes.
\end{theorem}

In this context the following definitions are convenient 
(for a motivation of these notions see \cite[Prop. 1.4]{Nil04}).

\begin{definition}{\rm
Let $P \subseteq \MR$ be a $d$-dimensional lattice polytope with $0 \in \intP$.
\begin{itemize}
\item $P$ is called a {\em canonical Fano} polytope, if $\intP \cap M = \{0\}$.
\item $P$ is called a {\em terminal Fano} polytope, if $P \cap M = \{0\} \cup \V(P)$.
\item $P$ is called a {\em smooth Fano} polytope, if the vertices of any facet of $P$ form a $\Z$-basis of 
the lattice $M$. 
\end{itemize}
}\end{definition}

If $P$ is a reflexive polytope, then $X_P$ is a nonsingular toric Fano variety if and only if $P^*$ is a 
smooth Fano polytope.

The following property \cite[Lemma 1.13]{Nil04} characterises reflexive polytopes.

\begin{lemma}
Let $P \subseteq \MR$ be a reflexive polytope.

For any $F \in \F(P)$ and $m \in F \cap M$ there is a $\Z$-basis $e_1, \ldots, e_{d-1}, e_d$ of $M$ such that 
$e_d = m$ and $F \subseteq \{x \in \MR \,:\, x_d = 1\}$; in particular 
$\eta_F = - e_d^*$ in the dual basis $e_1^*, \ldots, e_d^*$ of $N$. 

Furthermore $P$ is a canonical Fano polytope.
\label{root-local}
\end{lemma}

\section{The set of roots of a complete toric variety}

In this section the set of roots of a complete toric variety is investigated, 
and some classification results and bounds on the dimension of the automorphism group 
are achieved.

{\em Throughout the section let $\fan$ be a complete fan in $\NR$ with associated complete toric variety 
$X = X(N,\fan)$.}

\begin{definition}{\rm
Let $\J$ be the (finite) set of {\em Demazure roots} of $\fan$, i.e., 
\[\J := \{m \in M \;|\; \,\exists\, \tau \in \rays : \pro{\vt}{m} = -1 \text{ and } 
\,\forall\, \tau' \in \rays\backslash\{\tau\} :  \pro{\vts}{m} \geq 0\}.\]
For $m \in \J$ we denote by $\eta_m$ the unique primitive generator $\vt$ of a ray of $\fan$ 
with $\pro{\vt}{m} = -1$. 
For a subset $A \subseteq \J$ it is convenient to define $\eta(A) := \{\eta_m \,:\, m \in A\}$. 

Let $\S := \J \cap (-\J) = \{ m \in \J \::\: -m \in \J\}$ be the set of {\em semisimple} roots 
and $\calU := \J \backslash \S = \{m \in \J \::\: -m \not\in \J\}$ the set of {\em unipotent} roots. 
We say that $\fan$ is {\em semisimple}, if $\J = \S$, or equivalently $\calU = \emptyset$.

Furthermore we define $\S_1 := \{x \in \S \::\: \eta_x \not\in \eta(\calU)\}$ and $\S_2 := \S \backslash \S_1$, 
analogously $\calU_1 := \{x \in \calU \::\: \eta_x \not\in \eta(\S)\}$ and $\calU_2 := \calU\backslash \calU_1$. 
In particular $\eta(\S_1) \cap \eta(\S_2) = \emptyset$ and $\eta(\S_2) = \eta(\calU_2)$.
}
\label{root-defen}
\end{definition}

Usually the set $-\J$ is denoted as the set of Demazure roots (see \cite[Prop. 3.13]{Oda88}), 
however the sign convention here will turn out to be more convenient when considering normal fans of polytopes. 
Note that $\J$ only depends on the set of rays $\rays$.

For a root $m \in \J$ we can define a one-parameter subgroup 
$x_m \,:\, \C \to \Aut(X)$. Then the identity component $\Aut^\circ(X)$ is a semidirect product of 
a reductive algebraic subgroup containing the big torus $(\C^*)^d$ and having $\S$ as a root system and 
the unipotent radical that is generated by $\{x_m(\C) \,:\, m \in \calU\}$. Futhermore $\Aut(X)$ is generated by 
$\Aut^\circ(X)$ and the automorphisms that are induced by lattice automorphisms of the fan $\fan$. 
These results are due to Demazure (see \cite[p. 140]{Oda88}) in the nonsingular complete case, and were generalised 
by Cox \cite[Cor. 4.7]{Cox95} and B\"uhler \cite{Bue96}. 
Bruns and Gubeladze considered the case of a projective toric variety in \cite[Thm. 5.4]{BG99}. In particular 
there is the following result (recall that an algebraic group is reductive, if the unipotent radical is trivial).

\begin{proposition}
\begin{enumerate}
\item[]
\item $\Aut(X)$ is reductive if and only if $\fan$ is semisimple.
\item $\dim \Aut^\circ(X) = \card{\J} + d$.
\end{enumerate}
\label{autdim}
\end{proposition}

When $X_P$ is nonsingular, it is well-known (see \cite[p. 140]{Oda88}) that each irreducible component 
of the root system $\S$ is of type \textbf{A}. Here we will give an explicit description of $\S$ and $\eta(\J)$ 
by orthogonal families of roots that will turn out to be useful for applications.

\begin{definition}{\rm A pair of roots $v,w \in \J$ is called {\em orthogonal}, in symbols $v \bot w$, if 
$\pro{\eta_v}{w} = 0 = \pro{\eta_w}{v}$. In particular $\eta_{-v} \not= \eta_w \not= \eta_v \not= \eta_{-w}$.}
\label{ok}
\end{definition}

We remark that the term 'orthogonal' may be misleading, because most standard properties 
do not hold, e.g., $v \bot w$ does not necessarily imply $(-v) \bot w$.

\begin{lemma}
Let $B = \{b_1, \ldots, b_l\}$ be a non-empty set of roots such 
that\\$\pro{\eta_{b_i}}{b_j} = 0$ for $1 \leq j < i \leq l$. Then $B$ is a $\Z$-basis 
of $\lin_\R(B) \cap M$.
\label{ortho1}
\end{lemma}

\begin{proof}

We prove the base property by induction on $l$. Let $x := \lambda_1 b_1 + \cdots + \lambda_l b_l \in M$ 
with $\lambda_1, \ldots, \lambda_l \in \R$. Then $\lambda_l = - \pro{\eta_{b_l}}{x} \in \Z$. So $x - \lambda_l b_l = 
\lambda_1 b_1 + \cdots + \lambda_{l-1} b_{l-1} \in M$. Now proceed by induction.

\end{proof}

We define two special pairwise orthogonal families of roots:

\begin{definition}{\rm
Let $A \subseteq \J$. 

A pairwise orthogonal family $B \subseteq A$ is called 
\begin{itemize}
\item {\em $A$-facet basis}, if 
$\eta(A) = \{\eta_b \::\: b \in B\} \cup \{\eta_{-b} \::\: b \in B, -b \in A\}$.
\item {\em $A$-root basis}, if $A = \J \cap \lin(B)$.
\end{itemize}
\label{root-dfn}
}\end{definition}

When $B$ is an $A$-root basis, we have $\lin(A) = \lin(B)$, hence $\dim_\R \lin(A) = \card{B}$ by \ref{ortho1}. 
If furthermore $B \subseteq \S$, then Prop. \ref{root-fundi} below implies $A \subseteq \S$ and that $B$ 
is also an $A$-facet basis. 
Note that an $\S$-root basis is {\em not} a fundamental system for the root system $\S$ in the usual sense.

For an unambiguous description it is convenient to define an equivalence relation on the set of semisimple roots.

\begin{definition}{\rm
Let $v \equiv w$ ($v$ is {\em equivalent} to $w$), if $v,w \in \S$, $v \not= w$ and $\eta_{-v} = \eta_{-w}$.
In particular this yields $\pro{\eta_{-v}}{w} = -\pro{\eta_{-v}}{-w} = 1$.
\label{equidef}
}\end{definition} 

The goal of this section is to explicitly show how to get $\S$-root bases. 
To this end an algebraic-geometric approach due to Cox shall be discussed:

In \cite{Cox95} Cox described $\J$ as a set of {\em ordered} 
pairs of monomials in the homogeneous coordinate ring 
of the toric variety. 
For this we denote by $S := \C[x_\rho \,:\, \rho \in \rays]$ the homogeneous coordinate ring of $X$, i.e., 
$S$ is just a polynomial ring where any monomial in $S$ is naturally graded by 
the class group $\Cl(X)$, i.e., the degree of a monomial 
$\prod_\rho {x_\rho^{k_\rho}}$ is the class of the Weil divisor $\sum_\rho k_\rho \V_\rho$, where 
$\V_\rho$ is the torus-invariant prime divisor corresponding to the ray $\rho$. Recall that each 
$\rho \cap N$ is generated by $v_\rho$. 

We let $Y$ denote the set of indeterminates $\{x_\rho \,:\, \rho \in \rays\}$ and $\M$ the set of monomials in $S$. 
For any root $m \in \J$ we define $\rho_{m} := \pos(\eta_m) \in \rays$ and 
$x_m := x_{\rho_{m}} \in Y$. Now there is the following fundamental result 
\cite[Lemma 4.4]{Cox95} (with a different sign convention):

\begin{lemma}[Cox 95]
In this notation there is a well-defined bijection
\[h \;:\; \J \to \{(x_\rho,f) \in Y \times \M, \;:\; x_\rho \not= f, \, \deg(x_\rho) = \deg(f)\},\]
\[m \mapsto (x_m, \prod\limits_{\rho' \not= \rho_{m}} x_{\rho'}^{\pro{v_{\rho'}}{m}}).\]
For $m \in \J$ we have
\[m \in \S \lolra h(m) \in Y \times Y,\]
in this case $h(m) = (x_m, x_{-m})$.
\label{cox}
\end{lemma}

The next result can be used to 'orthogonalise' pairs of roots:

\begin{lemma}
Let $v,w \in \J$, $v \not= -w$, $\pro{\eta_v}{w} > 0$. Then 
$\pro{\eta_w}{v} = 0$, $v+w \in \J$ and $\eta_{v+w} = \eta_w$. 
Furthermore 
\[v+w \in \S \text{ iff } v \in \S \text{ and } w \in \S.\]
\label{tame}
\end{lemma}

\begin{proof}

Let $v$ correspond to $(x_v, f)$ and $w$ to $(x_w, g)$ as in Lemma \ref{cox}. It is $x_v \not= x_w$. 
The assumption implies that 
$x_v$ appears in the monomial $g$. Assume $\pro{\eta_w}{v} > 0$. Then $x_w$ would appear in the monomial $f$. 
However since $v \not= -w$ this is a contradiction to the antisymmetry of the order relation defined in 
\cite[Lemma 1.3]{Cox95}. The remaining statements are easy to see.

\end{proof}

\begin{corollary}
$v \in \calU$ and $w \in \S_1$ implies $\pro{\eta_v}{w} = 0$.
\label{neu}
\end{corollary}

Lemma \ref{tame} is a generalisation of parts of \cite[Prop. 3.3]{BG02} in a recent paper 
on polytopal linear groups due to Bruns and Gubeladze. 
The setting there is that of so called 'column structures' of polytopes where 'column vectors' correspond to roots. 
Most parts of this lemma were however already independently known and 
proven by the author as an application of Corollary \ref{root-prim-coro} below in the case of a reflexive polytope.

\begin{proposition}
Let $A \subseteq \J$ and $B \subseteq \S$ an $A$-root basis 
partitioned into $t$ equivalence classes of order $c_1, \ldots, c_t$. Then:

$\begin{array}{l}
A = \{\pm b \::\: b \in B\} \cup \{b - b' \::\: b,b' \in B,\, b \not= b',\, b \equiv b'\} \subseteq \S,\\
\card{A} = \card{B} + \sum_{i \in I} c_i^2 \leq \card{B} + \card{B}^2,\\
\eta(A) = \{\eta_{\pm b} \,:\, b \in B\},\; \card{\eta(A)} = \card{B} + t \leq 2 \card{B}.
\end{array}$

\label{root-fundi}
\end{proposition}

\begin{proof}

Only the first equation has to be proven: 
Let $m \in A$, by \ref{ortho1} we have 
$m = \sum_{b \in B} \lambda_b b$ for $\lambda_b \in \Z$. Let $l := \sum_{b \in B} \card{\lambda_b}$. 
Proceed by induction on $l$, let $l > 1$. 
By orthogonality we have $-1 \leq \pro{\eta_b}{m} = - \lambda_b$, hence $\lambda_b \leq 1$ for 
all $b \in B$. Assume there is an element $b \in B$ with 
$\lambda_b < 0$. Lemma \ref{tame} implies $b + m \in \lin(B) \cap \J = A$, so $b+m \in \S$ by induction hypothesis. 
Now Lemma \ref{tame} 
yields $-m \in A$. This yields $\lambda_b = -1$. Therefore $\lambda_b \in \{1,0,-1\}$ for all $b \in B$. 
Assume $l > 2$. By possibly replacing $m$ with $-m$ we can achieve that there are two elements $b,b' \in B$ with 
$\lambda_b = 1 = \lambda_{b'}$, hence $\eta_b = \eta_m = \eta_{b'}$, a contradiction. 
Therefore $l=2$, and there are two elements $b,b' \in B$ with $m = b - b'$. Assume $b \not\equiv b'$. Then 
necessarily $\pro{\eta_{-b'}}{b} = 0$, so $\eta_b = \eta_m = \eta_{-b'}$, a contradiction. 

\end{proof}

\begin{definition}{\rm
The grading of the polynomial ring $S := \C[x_\rho \,:\, \rho \in \rays]$ 
by the class group $\Cl(X)$ induces a partition of $Y$ into equivalence classes.
\begin{enumerate}
\item Let $Y_1, \ldots, Y_p$ be the equivalence classes of order at least two 
such that there exists no monomial in $\M\backslash Y$ of the same degree.
\item Let $Y_{p+1}, \ldots, Y_q$ be the remaining classes of order at least two.
\item Let $Y_{q+1}, \ldots, Y_r$ be the the equivalence classes of order one 
such that there exists an monomial in $\M\backslash Y$ of the same degree.
\item Let $Y_{r+1}, \ldots, Y_s$ be the remaining classes of order one.
\end{enumerate}}
\label{classes}
\end{definition}

By Lemma \ref{cox} ordered pairs of indeterminates contained 
in one of the equivalence classes $Y_1, \ldots, Y_p$ correspond to roots in 
$\S_1$, ordered pairs in $Y_{p+1}, \ldots, Y_q$ correspond to roots in $\S_2$. 
As changing $m \leftrightarrow -m$ for $m \in \S$ just means to reverse the corresponding pair of monomials, 
we immediately see that {\em no element in $\S_1$ is equivalent to an element in $\S_2$}. We have:
\[p = \card{\eta(\S_1)},\; q-p = \card{\eta(\S_2)}) = \card{\eta(\calU_2)},\; r-q = \card{\eta(\calU_1)},\; 
r = \card{\eta(\J)}.\]

We get from Lemma \ref{cox}:
\[\card{\S_1} = \sum\limits_{i=1}^p \card{Y_i} (\card{Y_i} - 1),\;\;\; 
\card{\S_2} = \sum\limits_{i=p+1}^q \card{Y_i} (\card{Y_i} - 1).\]

Moreover if we define for $i = p+1, \ldots, r$ the equivalence class $\M_i$ consisting of monomials in 
$\M\backslash Y$ having the same degree as an element in $Y_i$, we get:

\[\card{\calU_1} = \sum\limits_{i=q+1}^r \card{\M_i},\;\;\; 
\card{\calU_2} = \sum\limits_{i=p+1}^q \card{Y_i} \card{\M_i}.\]

In particular $\card{\calU_2} \not= \emptyset$ implies $\card{\calU_2} \geq 2$. 
Since by Lemma \ref{cox} for $i = p+1, \ldots, r$ no indeterminate in $Y_i$ can 
appear in an monomial in $\M_i$, we obtain that $v,w \in \calU$ with 
$\eta_v \not= \eta_w$ and $\deg(x_v) = \deg(x_w)$ are orthogonal. 
See Example \ref{weight} below for an illustration.

It is now simple to construct root bases:

\begin{proposition}
Let a subset $I \subseteq \{1, \ldots, q\}$ be given. Choose for any element $i \in I$ a 
subset $K_i \subseteq Y_i$ of cardinality $c_i + 1$. Denote by $R_i$ a set of $c_i$ semisimple roots 
corresponding to ordered pairs in $K_i$ with the same fixed second element. 
Define $B := \cup_{i \in I} R_i$ and $A := \lin(B) \cap \J$. 

Then $B$ is an $A$-root basis partitioned into equivalence classes $\{R_i\}_{i \in I}$, and any root in $A$ corresponds 
exactly to an ordered pair in $K_i$ for some $i \in I$.

Moreover any $A$-basis is given by this construction.
\label{root-constr}
\end{proposition}

\begin{proof}

By construction and Lemma \ref{cox} 
$\pro{\eta_v}{w} = 0 = \pro{\eta_w}{v}$ for $v,w \in B$, $v \not= w$, 
hence $B$ is an $A$-root basis with given equivalence classes. Using Lemma \ref{cox} and 
the description of $A$ in Prop. \ref{root-fundi} the remaining statements are easy to see.

\end{proof}

Choosing $I = \{1, \ldots, q\}$ and $K_i = Y_i$, so $c_i = \card{Y_i} - 1$ for $i \in I$, we get:

\begin{corollary}
$\S$-root bases exist, in particular $\J \cap \lin(\S) = \S$. Moreover
\[\dim_\R \lin(\S) = \sum\limits_{i = 1}^q (\card{Y_i} - 1).\]
\label{sort}
\end{corollary}

\begin{remark}{\rm 
It is interesting to note that $\conv(\S)$ is a centrally symmetric, terminal, 
reflexive polytope. More precisely due to \ref{root-fundi} and \cite[proof of Thm. 3.21]{DHZ01} 
there is an isomorphism of lattice polytopes (with respect to lattices $\lin(\S) \cap M$ and $\Z^{c_1 + \cdots + c_q}$)
\[\conv(\S) \cong (\Zt_{c_1} \oplus \cdots \oplus \Zt_{c_q})^*,\]
where $c_1, \ldots, c_q$ are defined as before, and 
$\Zt_n := \conv([0,1]^n, -[0,1]^n)$ is the $n$-dimensional standard zonotope. For a stronger statement 
see Theorem \ref{intersec}.
}\end{remark}

\begin{example}
{\rm 
Let's look at $X = \P^d$: 
We let $E_d$ denote the $d$-dimensional simplex $\conv(e_1, \ldots, e_d, -e_1 - \cdots - e_d)$, 
where $e_1, \ldots, e_d$ is a $\Z$-basis of $N$. Hence $E_d$ is 
the smooth Fano polytope corresponding to $d$-dimensional projective space $\P^d$. 
For $X = \P^d$ and $e^*_1, \ldots, e^*_d$ the dual basis of $M$ the family 
$b_1 := e^*_1,\, b_2 := e^*_1 - e^*_2,\, \ldots,\, b_d := e^*_1 - e^*_d$ forms an 
$\S$-root basis, where all elements are mutually equivalent. The homogeneous coordinate ring $S = \C[x_0, \ldots, x_n]$ 
is trivially graded. $\P^d$ is semisimple with $d^2+d$ roots.}
\label{pd-def}
\end{example}

\begin{example}{\rm 
For another example we consider the three-dimensional reflexive simplex 
$P := \conv((1,0,0),(1,3,0),(1,0,3),(-5,-6,-3))$ with $\V(P^*) = \{(-1,0,0), (-1,0,2), 
(2,-1,-1), (-1,1,0)\}$. We have $\dim_\R \S = 2$, $\card{\S} = 4$. 
$F_1$ and $F_2$ contain one antipodal pair of semisimple roots, while $F_3$ and $F_4$ contain the other pair. 
$F_3, F_4$ each contain three unipotent roots, pairs of unipotent roots in different facets 
are orthogonal. We can read this off the data 
$S = \C[x_0,x_1,x_2,x_3]$, $\Cl(X_P) \cong \Z$, $\deg(x_1) = \deg(x_2) = 1$ and 
$\deg(x_3) = \deg(x_4) = 2$. Hence $Y_1 = \{x_1,x_2\}$, $Y_2 = \{x_3,x_4\}$, $p=1$, 
$q=r=s=2$, $c_1 = 1 = c_2$, $\M_2 = \{x^2_1,x_1 x_2, x^2_2\}$. 
$X_P$ is just the weighted projective space with weights 
$(1,1,2,2)$.
\label{weight}}
\end{example}

Using above results we can show the existence of two special orthogonal families 
of roots (proof is left to the reader):

\begin{proposition}
\begin{enumerate}
\item[]
\item There exists an $\R$-linearly independent family $B$ of roots that can be partitioned into 
three pairwise disjoint subsets $B_1$, $B_2$, $B_3$ such that 
$B_1$ is an $\S_1$-root basis, $B_2$ is an $\S_2$-root basis, 
$B_1 \cup B_2$ is an $\S$-root basis and $B_3$ is a $\calU_1$-facet basis such that 
$\pro{\eta_b}{b'} = 0$ for all $b \in B_1 \cup B_2$ and $b' \in B_3$. 

Hence $\dim_\R \lin(\S) + \card{\eta(\calU_1)} = \card{B} \leq d$.
\item There exists an $\J$-facet basis $D$ that can be partitioned into 
three pairwise disjoint subsets $D_1$, $D_2$, $D_3$ such that 
$D_1$ is a $\calU_1$-facet basis, $D_2$ is a $\calU_2$-facet basis, 
$D_1 \cup D_2$ is a $\calU$-facet basis and $D_3$ is an $\S_1$-root basis. 

Hence $\card{\eta(\calU_1)} + \card{\eta(\calU_2)} + \dim_\R \lin(\S_1) = \card{D} \leq d$.
\end{enumerate}
\label{all}
\end{proposition}

There exists a classification result:

\begin{proposition}
A $d$-dimensional complete toric variety is isomorphic to a product of projective spaces iff 
there are $d$ linearly independent semisimple roots.

In this case
\[X \cong \P^{\card{Y_1}-1} \times \cdots \times \P^{\card{Y_q}-1}.\]
\label{charac}
\end{proposition}

\begin{proof}

Let $q=1$, so there is an $\S$-root basis $b_1, \ldots, b_d$ with $\eta_{-b_1} = \cdots = \eta_{-b_d}$. Assume 
there exists $\rho \in \rays$ with $\rho \not\in \{\rho_{b_1}, \ldots, \rho_{b_d}, \rho_{-b_1}\}$. Then 
$\pro{v_\rho}{b_i} = 0$ for $i = 1, \ldots, d$, since $b_i \in \S$. This implies $v_\rho = 0$, a contradiction. 
Therefore $\rays$ is determined. Since no cone in $\fan$ contains a linear subspace, this already implies 
$X \cong \P^d$. The general case is left to the reader. 

\end{proof}

\newpage
\begin{corollary}
\begin{enumerate}
\item[]
\item $\card{\eta(\calU)} \leq d$, $\card{\eta(\calU) \backslash \eta(\S)} \leq \codim_\R \lin(\S)$.
\item $\card{\eta(\J)} \leq 2 d$, with equality iff $X \cong \P^1 \times \cdots \times \P^1$.
\item $\card{\S} \leq d^2 + d$, with equality iff $X \cong \P^d$.
\end{enumerate}
\label{cories}
\end{corollary}

\begin{proof}

1. Follows from \ref{all}(2).

2. Let $D$ be the $\J$-facet basis from \ref{all}(2), we have $\card{D} \leq d$. 
By definition $\eta(\J) = \{\eta_{x} \::\: x \in D_1 \cup D_2\} \cup \{\eta_{\pm x} \::\: x \in D_3\}$, 
this gives the upper bound. Equality implies $D = D_3$, i.e., $\J = \S$, with 
no element in $D$ equivalent to any other. Applying the previous proposition we get the desired result.

3. Follows immediately from Corollary \ref{sort}, Prop. \ref{root-fundi} and Prop. \ref{charac}.

\end{proof}

While the case when $\MR$ is spanned by semisimple roots is completely classified, there are at least 
some partial results in the case of codimension one.

\begin{proposition}
Let $\dim_\R \lin(\S) = d-1$.

\begin{enumerate}
\item If $\card{\rays} \not= {\eta(\S)}$, then there exists $\tau \in \rays$ such that 
$\{\tau\} \subseteq \rays\backslash\eta(\S) \subseteq \{\pm \tau\}$, and we have 
$\Vt \cong \P^{\card{Y_1}-1} \times \cdots \times \P^{\card{Y_q}-1}$. 

\item If $q=1$, i.e., $\card{\S} = d^2 - d$, then $\card{\eta(\calU)} = 1$ and $\eta(\S) \cap \eta(\calU) = \emptyset$.
\end{enumerate}
\label{glocal}
\end{proposition}

\begin{proof}

Let $b_1, \ldots, b_{d-1}$ be an $\S$-root basis. By \ref{ortho1} we can find a lattice point $b_d \in M$ such 
that $b_1, \ldots, b_d$ is an $\Z$-basis of $M$. Let $e_1, \ldots, e_d$ denote the dual $\Z$-basis of $N$. 

1. Let $\tau \in \rays\backslash\eta(\S)$. Then $\pro{\vt}{b_i} = 0$ for all $i = 1, \ldots, d-1$, hence 
$\vt \in \{\pm e_d\}$. The set $\S$ is by construction canonically the set of roots of $\Vt$, so we can 
apply Prop. \ref{charac}.

2. Let $q=1$. By \ref{root-fundi} this is equivalent to $\card{\S} = (d-1)^2 + d-1 = d^2 - d$. 
For $i = 1, \ldots, d-1$ there exist $k_i \in \Z$ such that 
$\eta_i := 
\eta_{b_i} = - e_i + k_i e_d$. There exists $k_d \in \Z$ such that $\eta_d := 
\eta_{-b_1} = e_1 + \cdots + e_{d-1} + k_d e_d$. 

Since $\card{\eta(\S)} = d$, there exists $\tau \in \rays\backslash\eta(\S)$, we 
may assume $\vt = e_d$. 
Let $x = \lambda_1 b_1 + \cdots + \lambda_d b_d \in M$. We have $x \in \J$ with $\eta_x = e_d$ iff 
$\pro{x}{e_d} = \lambda_d = -1$ and $\pro{x}{\eta_i} \geq 0$ for $i = 1, \ldots, d$. This is equivalent to 
$\lambda_d = -1$, $\lambda_i \leq -k_i$ for $i = 1, \ldots, d-1$ and $\lambda_1 + \cdots + \lambda_{d-1} \geq k_d$. Hence 
there exists a root $x \in \J$ with $\eta_x = e_d$ if and only if $k_1 + \cdots + k_d \leq 0$.

On the other hand let $u := k_1 b_1 + \cdots + k_{d-1} b_{d-1} + b_d \in M$. Then $u^\bot$ is a hyperplane 
spanned by $\eta_1, \ldots, \eta_{d-1}$. We have $\pro{u}{e_d} = 1$ and $\pro{u}{\eta_d} = k_1 + \cdots + k_d$. 
Therefore if $\card{\rays} = d+1$, we get $\pro{u}{\eta_d} < 0$, so there exists 
$x \in \J$ with $\eta_x = e_d$, necessarily $e_d \in \eta(\calU)$. 
Otherwise if $\rays\backslash\eta(\S) = \{\pm e_d\}$, the analogous computation for $-e_d$ yields that 
either $e_d$ or $-e_d$ is in $\eta(\calU)$.

Assume $\eta(\S) \cap \eta(\calU) \not= \emptyset$, so $\S_2 \not= \emptyset$. 
Use the family $B$ in Prop. \ref{all}(1): Since by assumption 
all elements in $B_1 \cup B_2$ are mutually equivalent, however no element in $B_1$ 
is equivalent to one in $B_2$, we have $B = B_2$, i.e., $\S = \S_2$. This yields $\card{\eta(\calU_2)} = d$. Since 
$\card{\eta(\calU_1)} = 1$, we get a contradiction to \ref{cories}(1).

\end{proof}

For Gorenstein toric Fano varieties the second point cannot simply be improved as can be seen from Example 
\ref{weight}. Finally we obtain:

\begin{theorem}
Let $X$ be a $d$-dimensional complete toric variety with reductive automorphism group. Let 
$n := \dim \Aut^\circ(X)$. Then 
$$n \leq d^2 + 2 d,\; \text{ with equality only in the case of projective space.}$$

If $X$ is not a product of projective spaces, then 
$$n \leq d^2 - 2.$$
\label{root-dim}
\end{theorem}

\begin{proof}

Using \ref{autdim} we see that the first statement is just \ref{cories}(3). For the second statement use \ref{charac} 
to get $\dim_\R \S \leq d-1$, by \ref{root-fundi} we have in particular 
$\card{\S} \leq (d-1)^2 + d-1 = d^2 - d$. However equality cannot happen due 
to \ref{glocal}(2), since $\fan$ is semisimple. Since $\card{\S}$ is even, we get $\card{\S} \leq d^2 - d - 2$. 
Now use \ref{autdim}.

\end{proof}

\section{The set of roots of a reflexive polytope}

{\em Throughout the section let $P$ be a $d$-dimensional reflexive polytope in $\MR$.}

In this section we will focus on Gorenstein toric Fano varieties, these varieties correspond to reflexive polytopes 
as described in the first section. When $P$ is reflexive, we have by definition 
that {\em the set of roots $\J$ of the normal fan $\calN_P$ 
is exactly the set of lattice points in the relative interior 
of facets of $P$}.

\begin{definition}{\rm
The set $\J$ of roots of $P$ is defined as the set of roots of $\calN_P$. 
For $m \in \J$ we denote by $\F_m$ the unique facet of $P$ that contains $m$, 
and we again define $\eta_m = \eta_{\F_m}$ to be the unique primitive inner normal 
with $\pro{\eta_m}{\F_m} = -1$. 
For a subset $A \subseteq \J$ it is convenient to define $\F(A) := \{\F_m \,:\, m \in A\}$. 
We say $P$ is {\em semisimple}, if $\calN_P$ is semisimple, i.e., $\J = - \J$.
}\end{definition}

Most results of the previous section have now a direct geometric interpretation. Here three examples shall be 
explicitly stated.

\begin{corollary}
If a facet of a reflexive polytope contains an unipotent root and a semisimple root $x$, then 
the facet containing $-x$ also contains an unipotent root.
\end{corollary}

This follows from the fact that $\pm x$ corresponds to a pair of elements in one 
of the equivalence classes $Y_{p+1}, \ldots, Y_q$ (see \ref{classes}). Alternatively use \ref{neu}.

\begin{corollary}
Let $P$ be a $d$-dimensional reflexive polytope.

Then there are at most $2 d$ facets containing lattice points in their relative interior, equality holds 
iff $P \cong [-1,1]^d$ (isomorphic as lattice polytopes).
\label{wuerfel}
\end{corollary}

This follows from \ref{cories}(2). For another example we apply Prop. \ref{charac} and 
Prop. \ref{glocal}(2) to $d=2$ 
and get a characterisation of semisimple reflexive polygons without 
using the existing classification (e.g., \cite[Prop. 4.1]{Nil04}). 
The proof relies on the well-known fact that a two-dimensional terminal Fano polytope 
is a smooth Fano polytope, e.g., \cite[Lemma 1.17(1)]{Nil04}.

\begin{corollary}
Let $P$ be a two-dimensional reflexive polytope. For $k \in \N_{>0}$ let the reflexive polytope $E_k$ 
be defined as in \ref{pd-def}, i.e., $X_{E_k^*} \cong \P^k$. 

Then $P$ is semisimple iff $P$ is a smooth Fano polytope or $P \cong E_2^*$ or $P \cong E_1^2$.

$P$ or $P^*$ is semisimple iff $P$ or $P^*$ is a smooth Fano polytope.
\end{corollary}

To sharpen the results of the previous section we need an 
elementary but fundamental property of pairs of lattice points on the boundary of 
a reflexive polytope (for a proof see \cite[Prop. 3.1]{Nil04}).

\begin{lemma}
Let $v, w \in \randPM$. Then exactly one of the following is true:

\begin{enumerate} 
\item $v \sim w$, i.e., $v,w$ are contained in a common facet
\item $v + w = 0$
\item $v + w \in \randp$
\end{enumerate}

\noindent In the last case the following holds:

There exists exactly one pair $(a,b) \in \N^2_{>0}$ with 
$z := z(v,w) := a v + b w \in \randp$ such that $v \sim z$ and $w \sim z$. We have $a=1$ or $b=1$. 
Any facet containing $z$ (or $v+w$) contains exactly one of the points $v$ or $w$.

\label{root-prim}
\end{lemma}

The result shall be illustrated for $P := E_2^*$, i.e., $X_P \cong \P^2$:

\medskip

\centerline{\epsfxsize1.20in\epsffile{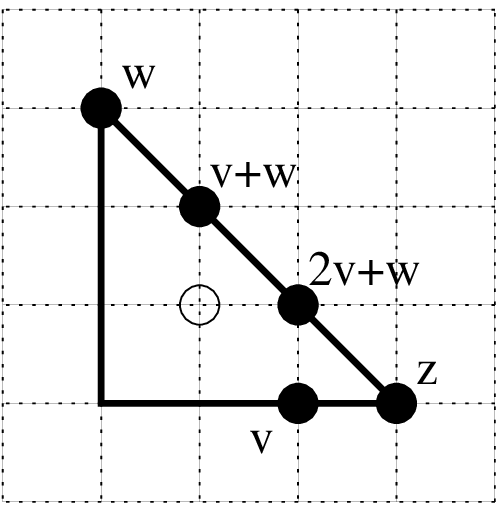}}

\medskip

This partial addition extends the partial 
addition of roots in \ref{tame} (see also \cite[Def. 3.2]{BG02}). Extending Definition \ref{ok} we may 
also more generally define $v \bot w$ for $v,w \in \randPM$, if $v + w \in \randp$ and $z(v,w) = v + w$.

\begin{corollary}
Let $v \in \J$, $w \in \randPM$ with $w \not\in \F_v$ and $w \not= -v$.

Then $v+w \in \randPM$ and $z(v,w) \in \F_v$. Moreover
$$\pro{\eta_v}{w} > 0 \text{ iff } z(v,w) = a v + w \text{ for } a \geq 2.$$
In this case $z(v,w) = (\pro{\eta_v}{w} + 1) v + w$.

\label{root-prim-coro}
\end{corollary}

There is a nice property of pairwise orthogonal families of roots:

\begin{proposition}
Let $B$ be a non-empty set of pairwise orthogonal roots. 

Then $F := \bigcap\limits_{b \in B} \F_b$ is a 
non-empty face of $P$ of codimension $\card{B}$, 
and the sum over all elements in $B$ is a lattice point in the relative interior of $F$.

\label{ortho2}
\end{proposition}

\begin{proof}

Let $B = \{b_1, \ldots, b_l\}$ with $\card{B} = l$. For $i \in \{1, \ldots, l\}$ we define 
$s_i := \sum_{j=1}^i b_j$ and $F_i := \cap_{j=1}^i \F_{b_j}$. 
Orthogonality implies that $\{\F_{b_1}, \ldots, \F_{b_l}\}$ is exactly the set of facets containing $s_l$. 
Therefore $s_l \in \relint F_l$, and since any $l$-codimensional face of $P$ is contained in at least 
$l$ facets, we must have $\codim F_l \leq l$. 
On the other hand $s_i \not\in F_{i+1}$ for all $i = 1, \ldots, l$, so $F_1 \supsetneq \cdots \supsetneq F_l$, 
hence we obtain $\codim F_l = l$.

\end{proof}

This yields a corollary concerning facets that contain unipotent roots. The proof follows from the existence of 
a $\calU$-facet basis (see \ref{all}(2)):

\begin{corollary}
If $\calU \not=\emptyset$, then $\bigcap\limits_{F \in \fb} F$ is a face of codimension $\card{\F(\calU)} \leq d$. 

In particular if $P$ is not semisimple, then 
the sum over all lattice points in the non-empty face $\bigcap_{F \in \fb} F$ is a non-zero fixpoint of $\Aut_M(P)$.
\label{ucoro}
\end{corollary}

Now we can improve Prop. \ref{charac} by taking the ambient space of semisimple roots into account 
(recall the definition of $E_d$ in \ref{pd-def}).

\begin{theorem}
Let $B \subseteq \S$ be an $A$-root basis for a subset $A \subseteq \J$, and 
$C_1, \ldots, C_t$ the partition of $B$ into equivalence classes of order $c_1, \ldots, c_t$. Then 
there are isomorphisms of lattice polytopes (with respect to lattices $\lin(A) \cap M$ and 
$\Z^{c_1 + \cdots + c_t}$)
$$P \cap \lin(A) \;\cong\; \bigoplus\limits_{i=1}^t P \cap \lin(C_i) \;\cong\; \bigoplus\limits_{i=1}^t 
E^*_{c_i}.$$
In particular the intersection of $P$ with the space spanned by all semisimple roots is again a reflexive polytope 
corresponding to a product of projective spaces.
\label{intersec}
\end{theorem}

\begin{proof}

Let $t=1$, i.e., all elements in $B$ are mutually equivalent. The general case is left to the reader. 
Let $l = \card{B} \geq 2$, $B = \{b_1, \ldots, b_l\}$, $b := b_1 + \cdots + b_l$.
$$\text{{\em Claim: }}\;\; P \cap \lin(b_1, \ldots, b_l) = \conv(b, b - (l+1) b_i \::\: i = 1, \ldots, l) \cong E_l^*.$$
Denote by $Q$ the simplex on the right hand side of the claim, so $Q \cong E_l^*$.

By \ref{ortho2} $b \in \bigcap_{i = 1}^l \F_{b_i}$. 
Since by assumption $\pro{\eta_{-b_i}}{b} = \sum_{j=1}^l \pro{\eta_{-b_i}}{b_j} = 
\sum_{j=1}^l \pro{\eta_{-b_j}}{b_j} = l$, it follows from \ref{root-prim-coro} that 
$z(-b_i, b) = b - (l+1) b_i \in \F_{-b_i}$ for $i = 1, \ldots, l$. Hence $Q \subseteq P \cap \lin(b_1, \ldots, b_l)$. 
On the other hand the previous calculation and orthogonality also implies 
that $Q \cap \F_{b_1}, \ldots, Q \cap \F_{b_l}, Q \cap \F_{-b_1}$ are exactly the facets of the simplex $Q$. 
This proves the claim. 

\end{proof}

\section{Criteria for a reductive automorphism group}

In this section we give several criteria for the automorphism group of a complete toric variety, respectively 
a Gorenstein toric Fano variety, to be reductive.

\begin{definition}{\rm For a polytope $Q \subseteq \MR$ we let $b_{Q}$ denote the {\em barycenter} of $Q$. 
When $Q$ is a lattice polytope, we denote by $\rvol(Q)$ 
the {\em relative lattice volume} of $Q$, i.e., $\rvol(\Pi) = 1$ for a fundamental paralleloped $\Pi$ 
of the lattice $\aff(Q) \cap M$.}
\end{definition}

\begin{theorem}
\begin{enumerate}
\item[]
\item Let $X = X(N,\fan)$ be a complete toric variety.

The following conditions are {\em equivalent}:
\begin{enumerate}
\item $\fan$ is semisimple, i.e., $\Aut(X)$ is reductive
\item $\sum\limits_{x \in \J} x = 0$
\item $\sum\limits_{\tau \in \rays} \pro{\vt}{x} = 0$ for all $x \in \J$
\end{enumerate}

If $\sum\limits_{\tau \in \rays} \vt = 0$, then $\fan$ is semisimple.
\item Let $X_P$ be a Gorenstein toric Fano variety for $P \subseteq \MR$ reflexive.

The following conditions are {\em equivalent}:
\begin{enumerate}
\item $P$ is semisimple, i.e., $\Aut(X_P)$ is reductive
\item $\sum\limits_{x \in \J} x = 0$
\item $\sum\limits_{v \in V(P^*)} \pro{v}{x} = 0$ for all $x \in \J$
\item $\sum\limits_{y \in P^* \cap N} \pro{y}{x} = 0$ for all $x \in \J$
\item $\pro{b_{P^*}}{x} = 0$ for all $x \in \J$
\item $\rvol(F') = \rvol(\F_x)$ for all $x \in \J$, $F' \in \F(P)$ with $\pro{\eta_{F'}}{x} > 0$
\item $\card{F' \cap M} = \card{\F_x \cap M}$ for all $x \in \J$, $F' \in \F(P)$ with $\pro{\eta_{F'}}{x} > 0$
\end{enumerate}

Any one of the following conditions is {\em sufficient} for $P$ to be semisimple:
\begin{enumerate}
\item[{\rm i.}] $b_{P} = 0$
\item[{\rm ii.}] $\sum\limits_{m \in P \cap M} m = 0$
\item[{\rm iii.}] $b_{P^*} = 0$
\item[{\rm iv.}] $\sum\limits_{y \in P^* \cap N} y = 0$
\item[{\rm v.}] $\sum\limits_{v \in V(P^*)} v = 0$
\item[{\rm vi.}] All facets of $P$ have the same relative lattice volume
\item[{\rm vii.}] All facets of $P$ have the same number of lattice points
\end{enumerate}

Condition {\rm vi} implies {\rm v}, e.g., if $P$ is a smooth Fano polytope.
\end{enumerate}
\label{main}
\end{theorem}

\begin{remark}{\rm Using the classification of $d$-dimensional reflexive polytopes for $d \leq 4$ due to 
Kreuzer and Skarke (see \cite{KS98,KS00,KS02}) we found examples showing that in the second part of the theorem 
the sufficient conditions 
${\rm i.} - {\rm v.}$ are pairwise independent, i.e., in general no condition implies any other.}
\end{remark}

\begin{remark}{\rm 
As explained in the introduction Batyrev and Selivanova obtained in \cite{BS99} from 
the existence of an Einstein-K\"ahler metric, that if $P$ is a reflexive polytope, $X_P$ is nonsingular 
and $P$ is {\em symmetric}, i.e., the group $\Aut_M(P)$ of linear lattice automorphisms leaving $P$ invariant 
has only the origin $0$ as a fixpoint, then $P$ has to be semisimple. They asked whether 
a direct proof for this combinatorial result exists. Indeed if $P$ is a symmetric reflexive polytope, we easily see 
that the second equivalent and even the first five sufficient conditions in the second part of the theorem 
are satisfied. For yet another proof see also Corollary \ref{ucoro}. Furthermore 
the first part of the theorem immediately yields a generalisation to complete toric varieties, for this see the fourth 
theorem in the introduction.}
\label{bs-sym}
\end{remark}

For the proof of Theorem \ref{main} we need some lemmas. The first is just a simple observation:

\begin{lemma}
Let $\fan$ be a complete fan in $\NR$.

\[m \in \J \lora \sum\limits_{\tau \in \rays} \pro{\vt}{m} \in \N,\]
in this case
\[m \in \S \lolra \sum\limits_{\tau \in \rays} \pro{\vt}{m} = 0.\]
\label{short}
\end{lemma}

\begin{lemma}
Let $\fan$ be a complete fan in $\NR$.

Let $A \subseteq \J$ be a subset such that 
\[\sum\limits_{m \in A} k_m m = 0\]
for some positive integers $\{k_m\}_{m \in A}$. Then $A \subseteq \S$.
\label{sum}
\end{lemma}

\begin{proof}

Assume $A \cap \calU \not= \emptyset$. Then by \ref{short}\\
$0 = \sum\limits_{\tau \in \rays} \pro{\vt}{\sum\limits_{m \in A} k_m m} = 
\sum\limits_{m \in A \cap \calU} k_m \sum\limits_{\tau \in \rays} \pro{\vt}{m} \geq 1,$ 
a contradiction.

\end{proof}

In the case of a reflexive polytope the following result is fundamental:

\begin{lemma}
Let $P$ be a $d$-dimensional reflexive polytope in $\MR$. 

Let $m \in \J$. Define the canonical projection map along m
\[\pi_m \,:\, M_\R \to M_\R / \R m.\]
Then $\pi_m$ induces an isomorphism of lattice polytopes
\[\F_m \to \pi_m(P),\]
with respect to the lattices $\aff(F) \cap M$ and $M/m \Z$.
\label{root-proj}
\end{lemma}

\begin{proof}

\cite[Prop. 2.2]{Nil04} immediately implies that $\pi_m :\, \F_m \to \pi_m(P)$ is 
a bijection. It is even an isomorphism of lattice polytopes by \ref{root-local}.

Another proof can be easily done using only the definition of a root.

\end{proof}

Using this lemma we get a reformulation of \ref{short}. Note that $A-B := \{a - b \,:\, a \in A, \,b \in B\}$ for 
arbitrary sets $A,B \subseteq \R^d$.

\begin{lemma}
Let $P$ be a $d$-dimensional reflexive polytope in $\MR$. 

For $m \in \J$ with $F := \F_m$ we have:

\begin{enumerate}

\item $P \subseteq F - \Rpos x$, $P \cap M \subseteq (F \cap M) - \N x$, 
$\{n \in P^* \cap N \,:\, \pro{n}{m} < 0\} = \{\eta_m\}$.

\item $P = \conv(F,F')$ iff there is only one facet $F'$ with $\pro{\eta_{F'}}{m} > 0$.

\item $m \in \S$ iff the previous condition is satisfied and $\pro{\eta_{F'}}{m} = 1$. 

In this case $F' = \F_{-m}$. Furthermore $F$ and $F'$ are naturally 
isomorphic as lattice polytopes and $\{n \in P^* \cap N \,:\, \pro{n}{m} \not= 0\} = \{\eta_m,\eta_{-m}\}$.
\end{enumerate}
\label{prisma}
\end{lemma}

\begin{lemma}
Let $P$ be a $d$-dimensional reflexive polytope in $\MR$. 

For $v \in \V(P^*)$ we denote by $v^* \in \F(P)$ the corresponding facet of $P$. Then
\[\sum\limits_{v \in \V(P^*)} \rvol(v^*)\, v = 0.\]
\label{mink}
\end{lemma}

\begin{proof}

Having chosen a fixed lattice basis of $M$ we denote by $\vol$ 
the associated differential-geometric volume in $\MR \cong \R^d$. Let $F \in \F(P)$ arbitrary. 
Since $\eta_F$ is primitive, it is a well-known fact that 
the determinant of the lattice $\aff(F) \cap M$, i.e., the volume of a 
fundamental paralleloped, is exactly $\norm{\eta_F}$, hence we get 
$\vol(F) = \rvol(F) \cdot \norm{\eta_F}$. The easy direction of the 
so called existence theorem of Minkowski (see \cite[no. 60]{BF71}) yields 
$\sum_{F \in \F(P)} \rvol(F)\, \eta_F = 0$.

\end{proof}

The approximation approach in the next proof is based upon an idea of Batyrev. 
Note that a facet $F$ of a $d$-dimensional 
polytope $Q \subseteq \MR$ is said to be parallel to $\R x$ for $x \in \MR$, if $\pro{\eta_F}{x} = 0$. 

\begin{lemma}
Let $Q \subseteq \MR$ be a $d$-dimensional polytope with a facet $F$ and $x \in \aff(F)$ such that 
$Q \subseteq F - \Rpos x$. For $q \in Q$ with $q = y - l x$, where $y \in F$ and $l \in \Rpos$, 
define $a(q) := y - 2 l x$. This definition extends uniquely to an affine map $a$ of $\MR$.

Then $a(b_{Q})$ is either in the interior of $Q$ or in the relative interior of a facet of $Q$ not parallel to $\R x$. 
The last case happens exactly iff there exists only one facet $F' \not= F$ not parallel to $\R x$.

\label{cone}
\end{lemma}

\begin{proof}

First assume there is exactly one facet $F' \not= F$ not parallel to $\R x$. This implies 
$Q = \conv(F,F')$. Choose an $\R$-basis $e_1, \ldots, e_d$ of $\MR$ such that $e_d = x$ and 
$\R e_1, \ldots, \R e_{d-1}$ are parallel to $F$. Now let $y \in F$ and define $h(y) \in \N$ such that 
$y - h(y) x \in F'$. For $k \in \Nspos$ let $F_k(y) := y + \cup_{i=1}^{d-1} \;[-\frac{1}{2 k}, \frac{1}{2 k}] e_i$ 
and $Q_k(y) := F_k(y) - [0,h(y)] x$. Then $b_{Q_k(y)} = y - \frac{h(y)}{2} x$ and 
$a(b_{Q_k(y)}) = y - h(y) x \in F'$. Let $M' := \Z e_1 + \cdots + \Z e_{d-1}$ and $z \in \relint F$. 
For any $k \in \Nspos$ we define $G_k := (z + \frac{1}{k} M') \cap F$ and $F_k := \cup_{y \in G_k} F_k(y)$. 
For $k \to \infty$ the sets $F_k$ converge uniformly to $F$. 
Therefore also $Q_k := \cup_{y \in G_k} Q_k(y)$ converges uniformly to $Q$ for 
$k \to \infty$. This implies that $b_{Q_k}$ converges to $b_{Q}$ for $k \to \infty$. 
Now $a$ is easily seen to be the restriction of an affine map of $\MR$, hence as 
$b_{Q_k}$ is a finite convex combination of $\{b_{Q_k(y)} \,:\, y \in G_k\}$ for any $k \in \Nspos$, also 
$a(b_{Q_k})$ is a finite convex combination of $\{a(b_{Q_k(y)} \,:\, y \in G_k\} \subseteq F'$ 
for any $k \in \Nspos$. This implies $a(b_{Q_k}) \in F'$ for any $k \in \Nspos$. Since $a$ is continuous and 
$F'$ is closed, this yields $a(b_{Q}) \in F'$. Hence obviously $a(b_{Q}) \in \relint F'$. 

Now let there be more than one facet different from $F$ that is not parallel to $\R x$. 
Then choose a polyhedral subdivision of $Q$ into finitely many polytopes $\{K_j\}$ such that any $K_j$ 
satifies the condition of the previous case. 
Then $b_{Q}$ is a proper convex combination of $\{b_{K_j}\}$, therefore by affinity 
also $a(b_{Q})$ is a proper convex combination of $\{a(b_{K_j})\} \subseteq \rand Q$. However since not all 
$a(b_{K_j})$ are contained in one facet, $a(b_{Q})$ is in the interior of $Q$. 

\end{proof}

\begin{proof}[Proof of Theorem \ref{main}]

The first part of the theorem, when $X$ is a complete toric variety, follows from \ref{short} and \ref{sum}. 
So let $X = X_P$ for $P \subseteq \MR$ a $d$-dimensional reflexive polytope, and we consider the second part of 
the theorem.

$(a)$ and $(b)$ are equivalent by \ref{sum}. 
The equivalences of $(a),(c),(d),(e)$ and the sufficiency of ${\rm iii,iv,v}$ follow from \ref{short} and 
\ref{prisma}.

$(f)$ and $(g)$ are necessary conditions for semisimplicity due to \ref{prisma}.

Let $(f)$ be satisfied and $x \in \J$. By \ref{prisma}(1) and \ref{mink} we have 
$$\rvol(\F_x) = \sum_{v \in \V(P^*),\, \pro{v}{x} > 0} \rvol(v^*) \pro{v}{x}.$$ 
By assumption there is only one 
vertex $v \in \V(P^*)$ with $\pro{v}{x} > 0$, furthermore $\pro{v}{x} = 1$. Hence \ref{prisma} implies $x \in \S$. 

Let $(g)$ be satisfied. Let $x \in \J$, $F := \F_x$ and $F' \in \F(P)$ with $\pro{\eta_{F'}}{x} > 0$. 
Due to \ref{prisma}(1) and by assumption 
there is a bijective map $h : F' \to F$ of lattice polytopes, i.e., $h(F' \cap M) 
\subseteq F \cap M$. Now there exists a lattice point $y \in F'$ with $h(y) = m$. 
Since $P = \conv(F,F')$ and $P$ is a canonical polytope, we obtain $y = -m \in \relint F'$, 
hence $m \in \S$.

The sufficiency of ${\rm vi,vii}$ is now trivial, \ref{mink} shows that ${\rm vi}$ implies ${\rm v}$. 

From now on let $x \in \J$ and $a$ the affine map defined as in \ref{cone} for $Q := P$ and $F := \F_x$. 

Let ${\rm i}$ be satisfied. 
By \ref{prisma}(1) we can apply Lemma \ref{cone} to get $-x = x-2x = a(0) = a(b_{P}) \in \J$, since 
$P$ is a canonical Fano polytope. 

Finally let ${\rm ii}$ be satisfied. For 
any $y \in F \cap M$ define $x_y \in P \cap M$ with $x_y := y - k x$ for $k \in \N$ maximal, and 
let $T_y := [y,x_y]$. Then \ref{prisma}(1) implies that

$\begin{array}{lcl}
-x \!\!\!&=&\!\!\! a(0) = a\left(\frac{1}{\card{P \cap M}} \sum\limits_{m \in P \cap M} m\right) = 
a\left(\sum\limits_{y \in F \cap M} 
\frac{\card{T_y \cap M}}{\card{P \cap M}} \frac{1}{\card{T_y \cap M}}\sum \limits_{m \in T_y \cap M} m\right)\\
&=&\!\!\! \sum\limits_{y \in F \cap M} 
\frac{\card{T_y \cap M}}{\card{P \cap M}} a\left(\frac{1}{\card{T_y \cap M}}\sum \limits_{m \in T_y \cap M} m\right) 
= \sum\limits_{y \in F \cap M} \frac{\card{T_y \cap M}}{\card{P \cap M}} x_y.
\end{array}$

Hence $-x$ is a proper convex combination of $\{x_y\}_{y \in F \cap M}$. 
As $P$ is a canonical Fano polytope, we get $-x \in \J$.

\end{proof}

\section{Centrally symmetric reflexive polytopes}

In this section the following result is going to be proved 
(recall the definition of the lattice polytope $E_1 := [-1,1]$).

\begin{theorem}
Let $P$ be a centrally symmetric $d$-dimensional reflexive polytope with $X_P$ the toric variety associated to $\calN_P$. 
Then
\begin{enumerate}
\item $P \cong E_1^{\frac{\card{\J}}{2}} \times G$ for a $\frac{\card{\J}}{2}$-codimensional 
face $G$ of $P$ that is a 
centrally symmetric reflexive polytope (with respect to $\aff(G) \cap M$ and a unique lattice point in 
$\relint\, G$) and has no roots itself.
\item Any facet contains at most $3^{d-1}$ lattice points and at most one root of $P$.\\
$P$ contains at most $3^d$ lattice points and has at most $2 d$ roots. Hence
\[\dim \Aut^\circ(X_P) \leq 3d.\]
\newpage
\item The following statements are equivalent:
\begin{enumerate}
\item $P$ contains $3^d$ lattice points
\item $P$ has $2 d$ roots, i.e., $\dim \Aut^\circ(X_P) = 3d$
\item Every facet of $P$ contains a root of $P$
\item Every facet of $P$ has $3^{d-1}$ lattice points
\item $P \cong E_1^d$, i.e., $X_P \cong \P^1 \times \cdots \times \P^1$
\end{enumerate}
\end{enumerate}
\label{central}
\end{theorem}

The first property immediately implies (see \ref{autdim}):

\begin{corollary}
Let $P$ be a centrally symmetric reflexive polytope with $X_P$ the toric variety associated to $\calN_P$. 

If $P$ contains no facet that is centrally symmetric with respect to a root of $P$, 
or there are at most $d-1$ facets of $P$ that 
can be decomposed as a product of lattice polytopes $E_1 \times F'$, then $P$ has no roots.

Hence if $d \geq 3$ and $P$ is simplicial, or $d \geq 4$ and any facet of $P$ is simplicial, 
then 
\[\dim \Aut^\circ(X_P) = d.\]
\end{corollary}

For the proof of Theorem \ref{central} we need the following lemma that 
is an easy corollary of \ref{prisma} and \ref{root-local}:

\begin{lemma}
Let $P$ be a centrally symmetric reflexive polytope.\\Let $F \in \F(P)$. Then
$$P \cong E_1 \times F \text{ iff } F \text{ contains a root } x \text{ of } P.$$
In this case $F$ is a centrally symmetric reflexive polytope (with respect to the lattice 
$\aff(F) \cap M$ with origin $x$). 
\label{rootcentral}
\end{lemma}

\begin{proof}[Proof of Theorem \ref{central}]

1. Apply the previous lemma inductively.

2. The bounds on the lattice points were proven in \cite[Thm. 6.4]{Nil04}. 
Since as just seen any facet of $P$ containing a root is reflexive, hence a canonical Fano polytope, 
it contains only one root of $P$. Now we apply \ref{autdim}(2) and 1. (or \ref{wuerfel}).

3. (b) $\lra$ (e) $\lra$ (c): Since $P$ as a centrally symmetric polytope contains at least $2 d$ facets, 
this follows from 1., alternatively use 2. and \ref{wuerfel}.

For the remaining equivalences we need the canonical map
$$\alpha \,:\, P \cap M \to M/3 M \cong (\GFt)^d.$$ 
In the proof of \cite[Thm. 6.4]{Nil04} it was shown that $\alpha$ is injective. 

Let $F \in \F(P)$ be arbitrary but fixed. Define $u := \eta_F \in \V(P^*)$ and also 
the $\GFt$-extended map $\alpha(u) : M/3 M \to \GFt$. For $m \in P \cap M$ 
we have $\pro{u}{m} \in \{-1,0,1\}$, in particular
\[m \in F \lolra \pro{\alpha(u)}{\alpha(m)} = -1 \in \GFt.\] 

(c) $\ra$ (a): Trivial, since (c) $\lra$ (e).

(a) $\ra$ (d): If $P$ contains $3^d$ lattice points, then $\alpha$ is a bijection, 
and therefore $\card{F \cap M} = \card{\{z \in M/3 M \,:\, \pro{\alpha(\eta_F)}{z} = -1\}} = 3^{d-1}$. 

(d) $\ra$ (c): The assumption implies that for any facet $F' \in \F(P)$ the map 
$$\alpha\res{F'} \,:\, F' \cap M \to 
\{z \in M/3 M \,:\, \pro{\alpha(\eta_{F'})}{z} = -1\}$$ 
is a bijection. Define 
$x := (1/{3^{d-1}}) \sum_{m \in F \cap M} m \in \relint F$. 

{\em It remains to prove $x \in M$.}

Choose a facet $G \in \F(P^*)$ and an $\R$-linearly independent family 
$w_1, \ldots, w_d$ of vertices of $G$ such that $w_1 = u$ and $w_2, \ldots, w_d$ are contained in a 
$(d-2)$-dimensional face of $P^*$.

Denote the corresponding facets of $P$ by $F_1, F_2, \ldots, F_d$ with $\eta_{F_j} = w_j$ for 
$j = 1, \ldots, d$, so $F_1 = F$. 
Then $Q := \cap_{j=2}^d F_j$ is a one-dimensional face of $P$. Therefore also the affine span 
of $\alpha(Q \cap M)$ is a one-dimensional affine subspace of $M/3 M$. Since $\card{F \cap Q} = 1$ 
there exists 
an element $b \in M/3 M$ such that $\pro{\alpha(u)}{b} = 0$ and $\pro{\alpha(w_j)}{b} = -1$ 
for all $j=2, \ldots, d$. 
Applying the assumption to $F_2$ yields a lattice point $v \in P \cap M$ with $\alpha(v) = b$. Hence also 
$\pro{u}{v} = 0$ and $\pro{w_j}{v} = -1$ for $j = 2, \ldots, d$. 

By \ref{root-local} we find a $\Z$-basis $e^*_1=u,e^*_2, \ldots, e^*_d$ of $N$ 
such that for any $j=2, \ldots, d$ there exist $\lambda_{j,k} \in \R$ with 
$e^*_j = \lambda_{j,2} (w_2-u) + \cdots + \lambda_{j,d} (w_d-u)$. 

\noindent$\cdot\;${\em Fact 1:} $\pro{w_k}{\sum_{m \in F \cap M} m} = 0$ for $k=2, \ldots, d$.

\noindent({\em Proof:} Since $F \cap F_k \not= \emptyset$, the assumption implies for $i=-1,0,1 \in \GFt$:
$\card{\{z \in M/3 M \,:\, \pro{\alpha(u)}{z} = -1, 
\pro{\alpha(w_k)}{z} = i\}} = 3^{d-2}$.)

\noindent$\cdot\;${\em Fact 2:} $\sum_{k=2}^d \lambda_{j,k} \in \Z$ for $j = 2, \ldots, d$.

\noindent({\em Proof:} 
$\pro{e^*_j}{v} $$=$$ (- \sum_{k=2}^d \lambda_{j,k}) \pro{u}{v} + \sum_{k=2}^d \lambda_{j,k} \pro{w_k}{v} $$=$ 
$ - \sum_{k=2}^d \lambda_{j,k}$ by the choice of $v$.)

Using these two facts we can finish the proof: 
\begin{eqnarray*}
\pro{e^*_1}{x} & = & \pro{u}{x} = -1 \in \Z,\\
\pro{e^*_j}{x} & = & (1/3^{d-1})\left((-\sum_{k=2}^d \lambda_{j,k}) \pro{u}{\sum_{m \in F \cap M} m} 
+ \sum_{k=2}^d \lambda_{j,k} \pro{w_k}{\sum_{m \in F \cap M} m}\right)\\ 
& = & \sum_{k=2}^d \lambda_{j,k} \in \Z \;\text{ for } j = 2, \ldots, d.
\end{eqnarray*}
Hence $x \in M$.

\end{proof}

\begin{remark}
{\rm Dropping the assumption of reflexivity and regarding just a complete toric variety $X = X(N,\fan)$ 
with centrally symmetric $\rays$ we still get immediately from \ref{root-defen}, \ref{autdim} and \ref{cories}(2) that 
$\dim \Aut^\circ(X) \leq 3d$, with equality if and only if $X \cong (\P^1)^d$.}

{\rm For $X$ as before, assume $X$ is also Gorenstein, i.e., the anticanonical divisor $-K_X$ is a Cartier divisor. In this 
case we can still show by slightly modifiying the proof of \cite[Thm. 6.4]{Nil04} that $h^0(X,-K_X) \leq 3^d$.}
\end{remark}

\end{document}